\documentclass[11pt]{article}

\usepackage[latin1]{inputenc}
\usepackage[french]{babel}
\usepackage{amsmath,amsthm,amsfonts,amssymb}
\usepackage{newlfont}
\usepackage{showidx}
\usepackage[dvips]{graphicx}
\usepackage{hyperref}

\newtheorem{theo}{Th\'eor\`eme}

\newtheorem{cor}{Corollaire}

\def\b{}

\def\g{\mathfrak g}
\def\p{\mathfrak p}
\def\k{\mathfrak k}
\def\n{\mathfrak n}

\def\b{\mathfrak b}
\def\s{\mathfrak s}

\def\hc{Harish-Chandra}
\def\ho{homomorphisme}

\def\ov{\overline}

\def\ad{\mathrm{ad}}

\begin{document}
\thispagestyle{empty}
\begin{center}
\Large \textbf{Paires symétriques orthogonales et Isomorphisme de
Rouvière}\vspace{1cm}

\normalsize  CHARLES TOROSSIAN

\vspace{1cm}

CNRS UMR 8553 , D.M.A.

École Normale Sup\'erieure,

 45 rue d'Ulm 75230 Paris
Cedex 05,

 torossia@dma.ens.fr

 10/09/2003
\end{center}

\noindent \textbf{Résumé:}  Dans un travail récent A. Alekseev et E.
Meinrenken (arXiv: math. RT/0308135) étudient les implications des
déformations de l'algèbre de Weyl en théorie de Lie. Ils montrent
notamment, que pour les paires symétriques orthogonales avec forme
bilinéaire anti-invariante pour l'involution de l'espace symétrique,
la formule  de Rouvière  est encore un isomorphisme d'algèbres. On
montre dans cette note comment retrouver ce résultat en utilisant
nos r\'esultats  sur l'homomorphisme d'Harish-Chandra généralisé
(\cite{Tor1}, \cite{Tor2} et \cite{Tor3}).

\section{Introduction}
Soit  $\mathfrak{g} $ une  algèbre de Lie  de dimension
 finie sur  $\mathbb{C}$, munie d'une involution $\sigma$ d'algèbre
 de Lie.
 On appellera $(\mathfrak{g},
\sigma)$ une paire  symétrique \cite{dix}. Si $\g$ n'est pas
algébrique au sens de Chevalley \cite{chev} alors on peut trouver
une paire symétrique $(\ov{\g}, \sigma)$ algébrique qui prolonge
$(\g, \sigma)$ comme paire symétrique avec   $\ov{\g}$ l'enveloppe
algébrique de $\g$ \cite{Tor4}.\\

 Soit
 $\mathfrak{g}=\mathfrak{k}\oplus \mathfrak{p} $  la décomposition en espaces propres relativement à
$\sigma$, on dira la  $\sigma$-décomposition. On notera
$\ov{\g}=\ov{\k}\oplus \ov{\p}$ la $\sigma$-décomposition de
$\ov{\g}$. Alors $\ov{\k}$ est l'enveloppe algébrique de $\k$
\cite{Tor4}. Soit $U(\mathfrak{g})$ l'algèbre enveloppante de
$\mathfrak{g}$. On définit pour  $X\in \mathfrak{k} $ le caractère
$\delta(X)=\frac{1}{2} {\mathrm tr}_{\mathfrak{g}/\mathfrak{k}}
{\mathrm{ad}}(X)$. On note $\mathfrak{k}^{-\delta}$ la sous-algèbre
de $U(\mathfrak{g})$ définie par
 $\{X-\delta(X), X\in \mathfrak{k}\}$.\\

 L'un des points fondamentaux en analyse harmonique est de
comprendre la nature de l'algèbre commutative (\cite{lich},
\cite{du79}) $(U(\mathfrak{g})/U(\mathfrak{g})\cdot
\mathfrak{k}^{-\delta})^{\mathfrak{k}}$
 , qui est directement liée à l'algèbre des opérateurs
différentiels invariants sur les demi-densit\'es de  l'espace
symétrique $G/K$ correspondant à ces données. \\

Nous avons défini dans \cite{Tor1}, ce que devrait être l'\ho\
d'\hc\ g\'en\'eralis\'e dans ce contexte général, de
$(U(\mathfrak{g})/U(\mathfrak{g})\cdot
\mathfrak{k}^{-\delta})^{\mathfrak{k} }$ sur ${\mathrm
S}[\mathfrak{p}]^{\mathfrak{k}}$ où $S[\mathfrak{p}]^{\mathfrak{k}}$
désigne l'espace des  $\mathfrak{k}$-invariants dans l'algèbre
symétrique $S[\mathfrak{p}]$. Pour l'instant notre formule est à
valeurs dans le corps de fractions rationnelles et nous conjecturons
(conjecture polynomiale) qu'elle est toujours à valeurs polynomiale
(voir aussi \cite{duflo-japon} pour des
conjectures plus générales).\\

Dans la plupart des cas connus (nous pensons que c'est toujours le
cas)  o\`u l'on sait que ces deux alg\`ebres sont isomorphes, notre
homomorphisme r\'ealise un isomorphisme d'alg\`ebres. En particulier
on retrouve l'isomorphisme de Rouvi\`ere du cas r\'esoluble et
l'isomorphisme de Helgason-\hc\
dans le cas des paires sym\'etriques r\'eductives.\\

Dans \cite{Tor3} nous avons \'etudi\'e le cas des paires
sym\'etriques orthogonales (c'est \`a dire munie d'une forme
bilin\'eaire invariante non d\'eg\'en\'er\'ee)  dont la forme
bilin\'eaire \'etait  $\sigma$-invariante. Nous avons montr\'e que
dans ce cas notre formule r\'ealisait encore un isomorphisme
d'alg\`ebres de $(U(\mathfrak{g})/U(\mathfrak{g})\cdot
\mathfrak{k}^{-\delta})^{\mathfrak{k} }$ sur ${\mathrm
S}[\mathfrak{p}]^{\mathfrak{k}}$, mais qu'en g\'en\'eral ce
n'\'etait pas
la formule de Rouvi\`ere.\\

En utilisant les résultats de notre th\`ese \cite{Tor1, Tor2} on
peut montrer (c.f. th\'eor\`eme \ref{2} de cette note) que la
formule de Rouvi\`ere s'\'etend aux paires sym\'etriques pour
lesquelles les \'el\'ements g\'en\'eriques de $\mathfrak{p}^*$
admettent des polarisations $\sigma$-stables v\'erifiant la
condition de Pukanszky. On notera $P-sigma$ cette condition.

Dans un travail récent A. Alekseev et E. Meinrenken \cite{AM}
étudient les implications des déformations de l'algèbre de Weyl en
théorie de Lie. Ils montrent notamment (th\'eor\`eme F), que pour
les paires symétriques orthogonales avec forme bilinéaire
anti-invariante pour l'involution $\sigma$, la formule de Rouvière
est encore un isomorphisme d'algèbres.\\

On montre dans cette note que la classe d'espaces sym\'etriques
consid\'er\'ee par A. Alekseev et E. Meinrenken, v\'erifie la
condition $ P-sigma$ (th\'eor\`eme \ref{1}) ce qui montre que la
formule de Rouvière s'applique dans ce cas.
\section{Une classe d'espaces sym\'etriques pour l'isomorphisme de Rouvière}
Rappelons dans un premier temps ce qu'est la formule de Rouvi\`ere.\\

Soit $X\in \p$, alors $\ad^2(X)$ est un endomorphisme de $\p$ et
l'on peut consid\'erer la fonction $J(X)$ d\'efinie par
\begin{equation}
J(X)=\det_\p\Big(\frac{\sinh \ad X}{\ad X}\Big).
\end{equation}

La s\'erie formelle  $J^{1/2}$ est bien d\'efinie car $J(0)=1$,
c'est un \'el\'ement de $S[[\p^*]]$.  On peut alors consid\'erer
$\partial_{J^{1/2}}$ l'op\'erateur diff\'erentiel d'ordre infini \`a
c\oe fficients constants agissant comme endomorphisme sur $S[\p]$.\\

On note $\beta$ la sym\'etrisation de $S[\p]$ sur $U(\g)/U(\g)\cdot
\k^{-\delta}$ . Pour $P\in S[\p]^{\k}$ on note \cite{rou86}
$$R(P)=\beta(\partial_{J^{1/2}}P).$$ C'est un \'el\'ement de
$(U(\g)/U(\g)\cdot \k^{-\delta})^{\k}$ et $R$ est  la formule de
Rouvi\`ere. C'est une  formule analogue \`a celle de Duflo
\cite{du77}. Cette formule d\'efinit dans le cas des paires
sym\'etriques r\'esolubles \cite{rou86} ou dans le cas des paires
sym\'etriques v\'erifiant la condition $P-sigma$, un isomorphisme
d'alg\`ebres de $S[\p]^{\k}$ sur
$(U(\g)/U(\g)\cdot \k^{-\delta})^{\k}$ (th\'eor\`eme \ref{2}).\\

\subsection{Construction de polarisations $\sigma$-stables}

Le premier r\'esultat de cette note est:

\begin{theo}\label{1} Soit $(\g, \sigma)$ une paire sym\'etrique orthogonale
sur $\mathbb{C}$, munie d'une
 forme bilin\'eaire invariante, non d\'eg\'en\'er\'ee et
anti-invariante pour $\sigma$. On note $(\ov{\g}, \sigma)$ une
enveloppe algébrique comme paire symétrique, alors  tout \'el\'ement
g\'en\'erique   $f\in \ov{\p}^*$ admet une polarisation
$\sigma$-stable v\'erifiant la condition de Pukanszky.
\end{theo}

Avant de d\'emontrer le th\'eor\`eme faisons quelques rappels sur
les \'el\'ements g\'en\'eriques de $\p^*$ dans notre situation.

On note $B$ la forme bilinéaire invariante et non dégénérée. Par
d\'efinition de l'anti-invariance pour $\sigma$ on a  $$B(\sigma(x),
\sigma(y))=-B(x,y).$$ On en d\'eduit imm\'ediatement que $\k$ et
$\p$ sont isotropes et en dualit\'e. En d\'efinitive $\p^*$ est
isomorphe comme $\k$-module \`a $\k$.\\

Soit $f\in \p^*=\k^\perp$ une forme lin\'eaire nulle sur $\k$.
Notons $x_f$ l'\'el\'ement de $\k$ correspondant via la forme $B$.
Cet \'el\'ement admet une d\'ecomposition de Jordan dans $\ov{\k}$,
notons l\`a
$$x_f=x_s+x_u.$$Le noyau de la forme altern\'ee $B_f$ (rappelons que
l'on a $B_f(x,y)=f[x,y]$) associ\'ee \`a $f$ est alors le
centralisateur de $x_f$, on le note comme d'habitude $\g^{f}$. Alors
$\g^f$ est $\sigma$-stable et on a $\g^f=\k^f\oplus\p^f$ . De plus
$\k/\k^f$ et $\p/\p^f$ sont en dualité par $B_f$ \cite{Kost-Rallis}.
Comme $\k$ et $\p$ ont même dimension
on en déduit que c'est aussi le cas pour $\k^f$ et $\p^f$.\\

Un \'el\'ement $f\in \ov{\p}^*$ est r\'egulier au sens des paires
symétriques si $\ov{\g}^f$ est de dimension minimale parmi les
$\ov{\g}^g$ pour $g\in  \ov{\p}^*$. Les formes lin\'eaires
r\'eguli\`eres dans $\ov{\p}^*$ forment un ouvert de Zariski et
v\'erifient les conditions suivantes (\cite{Tor1} paragraphes 1.2 et
1.3):
 \begin{enumerate}
 \item $[\ov{\k}^f, \ov{\p}^f]=0$.
 \item la sous-algèbre engendrée par $\ov{\p}^f$ est nilpotente et on
 note $\s_f$ son tore maximale. Il est inclus dans $\ov{\p}$.

 \end{enumerate}

On dit que $f$ est  g\'en\'erique (on dit aussi tr\`es r\'egulier)
si le tore $\s_f$ est de dimension maximale parmi les tores $\s_g$
avec $g$ régulier. Les tores $\s_f$ pour $f$ g\'en\'eriques sont
tous conjugu\'es (\cite{Tor1} proposition
1.4.2.1).\\

Une polarisation en $f$ dans $\ov{\g}$ est une sous-algèbre isotrope
pour $B_f$ et de dimension  maximale. C'est automatiquement une
sous-algèbre algébrique. On dit que $f$ vérifie la condition de
Pukanszky si on a $\b=\ov{\g}^f+ \b_u$ avec $\b_u$ le radical
unipotent de $\b$.

Notons $Z(\ov{\g})$ le centre de $\ov{\g}$ et $Z(\ov{\g})_s$ sa
partie semi-simple. Le $\sigma$-rang de $\ov{\g}$ est d\'efini comme
\'etant la dimension de
tores $\s_f/\s_f\cap Z(\ov{\g})_s$  pour $f$ g\'en\'erique dans $\ov{\p}^*$.\\

\noindent \textbf{Preuve du th\'eor\`eme \ref{1}:} La
d\'emonstration se fait comme
d'habitude par r\'ecurrence sur la dimension de $\g$.\\

La forme $B$ est invariante par $\g$ et donc  par $\ov{\g}$. Soit
$f\in \ov{\p}^*$, g\'en\'erique. On note encore $f$ sa restriction à
$\g$. On \'ecrit $x_f=x_s+x_u$. On voit facilement que $\ov{\g}^f$
est encore le centralisateur de $x_f$ dans $\ov{\g}$. En effet on a
pour tout $x \in \ov{\g}$ ,  $$[x, \ov{\g}]=[x, \g].$$Alors  si
$f([x, \ov{\g}])=f([x, \g])=0$, on a $B(x_f, [x, \g])=B([x_f,
x], \g)=0$, ce qui permet de conclure.\\

Consid\'erons la d\'ecomposition de $\g$ sous l'action adjointe de
$x_s$ (rappelons que $\g$ est un idéal de $\ov{\g}$), On a
\begin{equation}\g=\g_o\oplus \sum_{\lambda}
\g_{\lambda}\end{equation} avec $\g_{\lambda}$ espace propre pour la
valeur propre $\lambda$. Comme on a $\sigma(x_s)=x_s$ l'espace
propre $\g_{\lambda}$ est $\sigma$-stable. Comme on a
$\ov{\g}_{\lambda}=\g_{\lambda}$ on a aussi
\begin{equation}\label{decomposition}\ov{\g}=\ov{\g}_o\oplus \sum_{\lambda}
\g_{\lambda}.\end{equation}On a clairement $\ov{\g}_o=\ov{\g_o}$ enveloppe algébrique de $\g_o$.\\

La forme $B$ \'etant invariante par $\ov{\g}$, on a
$\g_{\lambda}\perp_{B} \g_{\mu}$ si $\lambda +\mu \neq 0$. On en
d\'eduit que si $\lambda$ est racine alors $-\lambda$ l'est aussi et
que $\g_{\lambda}$ et $\g_{-\lambda}$ sont en dualit\'e. Comme dans
\cite{Tor1} (comme dans le cas réductif) on peut s\'eparer les
racines en deux sous-ensembles $\Delta$ et $-\Delta$ stables par
addition. Alors $\mathfrak{n}=\sum_{\lambda \in \Delta} \g_\lambda$
est une
sous-alg\`ebre $\sigma$-stable.\\

La sous-alg\`ebre $\g_o$ est $\sigma$-stable, contient $\g^f$ et
$B|{\g_o\times \g_o}$ est non d\'eg\'en\'er\'ee. Par conséquent
$(\g_o, \sigma)$ est un paire symétrique qui est dans la même classe
que $\g$. Pour construire une polarisation $\sigma$-stable de $f$
dans $\ov{\g}$ il suffit de construire une polarisation
$\sigma$-stable de $f$ dans $\ov{\g}_o$ vérifiant la condition de
Pukanszky et de lui ajouter
$\n$.\\

La décomposition (\ref{decomposition}) est orthogonale pour $B_f$ et
$B_f$ est non dégénérée sur $[x_s, \ov{\g}]$. On en déduit
facilement que $f|_{\ov{\g}_o}$ est encore g\'en\'erique dans
$\ov{\p}_o^*$ (\cite{Tor1}, paragraphe 1.4). On pourra donc
appliquer l'hypoth\`ese de r\'ecurrence si la dimension de
$\ov{\g}_o$ est inf\'erieure \`a la dimension de $\ov{\g}$, c'est à
dire si
$\ad_{\g}(x_s)$ est non nul.\\

Supposons que pour $f$ générique on ait  $\ad_{\g}(x_s)=0$. Je dis
alors que $\ov{\g}$ est de $\sigma$-rang nul et d'après le
th\'eor\`eme 1.5.3.1 de \cite{Tor1} cet élément admet une
polarisation $\sigma$-stable vérifiant la condition de Pukanszky.

En effet si le tore $\s_f$ n'est pas central dans $\ov{\g}$ on a une
d\'ecomposition en espaces radiciels sous l'action adjointe de
$\s_f$
$$\ov{\g}=\ov{\g}_{(o)}\oplus \sum_{\mu} \g_{(\mu)}.$$
Comme $f$ est régulier on a  $\ov{\g}^f\subset \ov{\g}_{(o)}$. Alors
$\g_{(\mu)}$ est stable par $\ov{\g}_{(o)}$, en particulier sous
l'action adjointe $\ad(x_f)$ car $x_f\in \ov{\g}^f\subset\g_o$. Or
$\ad(x_f)$ est unipotent (car $\ad(x_s)=0$) donc il existe $y\in
\g_{\mu}$ non nul v\'erifiant $[x_f, y]=0$, c'est à dire $y\in \g^f
\subset \ov{\g}_{(o)}$ ce qui est absurde, car $\g_{(o)}\cap
\g_{\mu}$ est
nul. $\Box$\\

\noindent \textbf{Remarque:} Si $(\ov{G},\sigma)$ est un groupe
alg\'ebrique connexe correspondant \`a la paire sym\'etrique
$(\ov{\g}, \sigma)$, notre d\'emonstration montre facilement que
l'on peut choisir une polarisation stable par $\ov{K}(f)$ avec
$\ov{K}$ les points de
$\sigma$ dans $\ov{G}$.\\

\subsection{Une extension des r\'esultats de \cite{Tor1, Tor2}}

Nous consid\'erons que le deuxi\`eme r\'esultat de cette note
\'etait acquis et contenu dans \cite{Tor1, Tor2}. Par souci de
claret\'e nous pr\'esentons ici ce r\'esultat sous une forme plus
autonome:\\

\begin{theo}\label{2} Soit $(\g, \sigma)$ une paire symétrique
alg\'ebrique sur $\mathbb{C}$ v\'erifiant la condition $P-sigma$,
alors la formule de Rouvi\`ere d\'efinit un isomorphisme
d'alg\`e\-bres de $S[\p]^{\k}$ sur $(U(\g)/U(\g)\cdot
\k^{-\delta})^{\k}$.
\end{theo}

Le principe de l'extension des scalaires assure que ce résultat
est encore vrai sur le corps des r\'eels. Comme $\g$ est un idéal
dans $\ov{\g}$, on d\'eduit de ces  deux th\'eor\`émets
pr\'ec\'edents
une autre preuve du th\'eor\`eme F de \cite{AM}. \\

\begin{cor}[\cite{AM}]
Soit $(\g, \sigma)$ une paire sym\'etrique orthogonale avec forme
bilin\'eaire anti-invariante par $\sigma$, alors la formule de
Rouvi\`ere est un isomorphisme d'alg\`ebres de $S[\p]^{\k}$ sur
$(U(\g)/U(\g)\cdot \k^{-\delta})^{\k}$.
\end{cor}

\noindent \textbf{Preuve du th\'eor\`eme \ref{2} :} Ce th\'eor\`eme
est une conséquence non triviale de la méthode des orbites et c'est
l'essence de la preuve de Duflo \cite{du77} que l'on a appliquée aux
espaces symétriques dans \cite{Tor1, Tor2} (voir aussi \cite{ben} et
\cite{cimpa}).

On ébauche les grandes lignes de la démonstration qui consiste en
une redite des textes \cite{Tor1, Tor2}.

Introduisons quelques notations indispensables dont on trouvera les
d\'efini\-tions d\'etaill\'ees dans \cite{Tor1, Tor2}.

 On se donne
$(G, \sigma)$ un groupe lin\'eaire alg\'ebrique muni d'une
involution vérifiant  $Lie(G)=\g$ et tel que la diff\'erentielle de
$\sigma$ \`a l'origine co\" \i ncide avec l'involution de
l'alg\`ebre de Lie (que l'on note par la m\^eme lettre). On pose
$K=G^{\sigma}$ le sous-groupe des points fixes. On note $\Delta_G$
le caract\`ere $|\det_{\g}Ad (x)|^{-1}$ et
$\Delta_{G,K}(x)=\Delta_K(x)/\Delta_G(x)=|\det_{\g/\k}Ad(x)|$ le
caract\`ere de $K$.

 Pour $\chi$ un caract\`ere de $K$ on note
$G\times_K \mathbb{C}_{\chi}$ le fibr\'e en droites dont les
sections v\'erifient $\phi(gk)=\chi(k)^{-1}\phi(g)$. On d\'efinit
(voir par exemple \cite{Tor1}, paragraphe 2.1) via les mesures de
Haar une densit\'e \`a valeur dans $G\times_K
\mathbb{C}_{\Delta_{G,K}}$ not\'e $d_{G/K}$ . On a alors
$$d_{G,K}(gk)=\Delta_{G,K}(k)^{-1}d_{G,K}(g).$$
On considère les objets réels sous-jacents. Comme $G$ est un
sous-groupe d'un groupe lin\'eaire on  introduit l'ouvert invariant
et $\sigma$-stable, $\mathcal{V}$ des \'el\'ements de $X\in \g$
v\'erifiant $|Im(\lambda)|< \pi/4$ pour toutes valeurs propres
$\lambda$ de $X$. L'application exponentielle est un difféomorphisme
de $\mathcal{V}\cap \p $ sur son image dans  $G/K$, de plus
$\mathcal{W}=\p\cap \mathcal{V}$ est
invariant par l'action adjointe de $K$.\\

La mesure de Lebesgue sur $\p^*$ est semi-invariante par $K$ de
poids $\Delta_{G,K}$. Par d\'esint\'egratrion on obtient sur presque
toutes les orbites g\'en\'eriques $\omega$ une mesure de  poids
$\Delta_{G,K}$ que l'on  note $d\lambda_{\omega}$. On peut sans
perte de g\'en\'eralit\'e supposer que ces mesures sont
temp\'er\'ees, ce qui permet de consid\'erer leurs transform\'ees de
Fourier que l'on multiplie par $J^{-1/2}$. Ce sont des fonctions
g\'en\'eralis\'ees sur $\mathcal{W}$ semi-invariantes
de poids $\Delta_{G,K}^{-1}$.\\

 Le point clef est de montrer que ces
fonctions généralisées sont  propres sous l'action naturelle de
$(U(\g)/U(\g)\cdot \k^{-\delta})^{\k}$, alg\`ebre
des op\'erateurs diff\'erentiels invariants sur les demi-densit\'es de $G/K$.\\

On consid\`ere l'application exponentielle not\'ee $Exp$ de $\p$ sur
$G/K$. Comme c'est un difféomorphisme de $\mathcal{W}$ sur son image
et  on peut consid\'erer l'image r\'eciproque du fibr\'e en droites
$Exp_\mathcal{W}^*(G\times_{K}\mathbb{C}_{\Delta_{G,K}^{1/2}})=\mathcal{W}\times
\mathbb{C}$. \\

Soit $\varphi(X)dX$ une densit\'e sur $\mathcal{W}$ \`a support
compact et consid\'erons  la densit\'e du fibr\'e
$G\times_{K}\mathbb{C}_{\Delta_{G,K}^{1/2}}$ \`a support dans
$Exp(\mathcal{W})$ d\'efinie par
$$Exp_*(\varphi(X)dX)=\Phi(g)d_{G,K},$$ alors $\Phi$ est une section
du fibr\'e $G\times_{K}\mathbb{C}_{\Delta_{G,K}^{-1/2}}$.\\

Par hypoth\`ese les orbites génériques admettent des polarisations
$\sigma$-stables qui vérifient la condition de Pukanszky. Pour un
tel $f\in \p^*$ on note $\b$ une polarisation $\sigma$-stable.
Notons $B$ le groupe alg\'ebrique connexe d'alg\`ebre de Lie $\b$.
Notons $\chi_f$ la fonction d\'efinie sur $B\cap
\exp(\mathcal{V})=\exp(\b\cap \mathcal{V})$ (\cite{Tor1}, lemme
2.3.3.1), par $\chi_f(exp(X))=e^{if(X)}$. Remarquons que l'on n'a
pas besoin de d\'efinir un caract\`ere du groupe $B$ ce qui a
l'avantage de ne pas supposer que $if|_{\b}$ soit int\'egrable. On
consid\`ere la section g\'en\'eralis\'ee du fibr\'e
$G\times_{K}\mathbb{C}_{\Delta_{G,K}^{-1/2}}$ d\'efinie sur
$Exp(\mathcal{W})$  par
\begin{equation}\Phi d_{G,K}\longmapsto j_*(\Phi d_{G,K})=\int_{B/B\cap K}
\Phi(b)\chi_f(b)\Delta_{G,B}^{-1/2}(b) d_{B, B\cap K}(b).
\end{equation} D'apr\`es \cite{Tor1}, 2.8 tous les caract\`eres sont
bien ajust\'es. D'apr\`es \cite{Tor1}, th\'eor\`eme 3.3.2.1, c'est
une section g\'en\'eralis\'ee propre sous l'action des op\'erateurs
diff\'erentiels invariants $D_u$ pour $u\in (U(\g)/U(\g)\cdot
\k^{-\delta})^{\k}$. Ceci fournit donc un
caract\`ere de cette alg\`ebre not\'e $u \longrightarrow\lambda_{f,\b}(u)$.\\

Pour $k\in K$  la formule
$$\Phi d_{G,K}\longmapsto j_*(\Phi d_{G,K})(k)=\int_{B/B\cap K}
\Phi(kb)\chi_f(b)\Delta_{G,B}^{-1/2}(b) d_{B, B\cap K}(b)$$d\'efinit
encore une section g\'en\'eralis\'ee propre sous l'action des
op\'erateurs $D_u$ de m\^eme caract\`ere $\lambda_{f,\b}$ car $D_u$
est invariant.

La condition de Pukanszky assure que l'orbite $K\cdot f$ est
fibr\'ee
$$K\cdot
f=K\times_{K\cap B} (f+(\p\cap \b)^{\perp}).$$ Les paragraphes 2.6
et 2.7 de \cite{Tor1} montrent alors que
\begin{equation}
\int_{K/K\cap B} j_*(\Phi d_{G,K})(k) \Delta_{G,K}^{-1/2}(k)
d_{K,K\cap B}(k)
\end{equation} est proportionnelle \`a
\begin{equation}
\int_{K\cdot f} d\lambda_{\omega}(\xi)\left(\int_{\p} \varphi(X)
e^{i\xi(X)}J^{-1/2}(X)\right).
\end{equation}
On conclut d'une part que la transform\'ee de Fourier de l'orbite
$\omega=K\cdot f$ multipli\'ee par $J^{-1/2}$ est une fonction
g\'en\'eralis\'ee sur $\mathcal{W}$ propre sous l'action des
op\'erateurs $D_u$ écrits en coordonn\'ees exponentielle. Comme on
peut reconstituer la masse de Dirac en $0\in \p$, on en d\'eduit
comme dans \cite{du77, ben, Tor2} que pour presque tout $f\in \p^*$
le caract\`ere $\lambda_{f,\b}$ est donn\'e par la formule de
Rouvi\`ere. Ceci montre que la formule de Rouvi\`ere est un
homomorphisme d'alg\`ebres.$\Box$.

\subsection{Comparaison avec l'homomorphisme d'Harish-Chandra
g\'en\'eralis\'e}

Expliquons maintenant pourquoi dans l'hypoth\`ese d'une paire
sym\'etrique orthogonale avec forme bilin\'eaire anti-invariante par
$\sigma$, notre homomorphisme d'\hc\ g\'en\'eralis\'e co\" \i
nciderait  avec la formule de Rouvi\`ere.\\

La démonstration de  \cite{Tor3} concernant le cas des paires
sym\'etriques r\'esolu\-bles est bas\'ee sur l'\'etude des parties
radiales des op\'erateurs diff\'erentiels invariantes "\`a
l'infini". C'est ce que l'on fait dans le cas des espaces
sym\'etriques r\'eductifs.  La comparaison des facteurs dominants
fournit dans le
cas r\'esoluble le r\'esultat cherch\'e.\\

Cette d\'emonstration  fonctionne sans problème pour le cas
\'etudi\'e dans cette note sous la condition que l'on sache écrire
l'action des opérateurs différentiels invariants en coordonnées
exponentielles sur les distributions $K$-semi-invariantes. A savoir
après modification par $J^{1/2}$, ces opérateurs agissent comme des
opérateurs à c\oe fficients constants, ce qui expliquerait a
posteriori que les transformées de Fourier des $K$-orbites dans
$\p^*$ soient des
fonctions propres comme nous venons de le voir.\\

C'est le  problème de l'extension de l'isomorphisme de Rouvière (ou
celui de Duflo dans le cas des groupes)  au cas des distributions à
support quelconque, qui motive les conjectures de
Kashiwara-Vergne-Rouvière sur les formules de
Baker-Campbell-Hausdorff (\cite{KV, rou90, rou91, rou94, ADS, AST,
Tor5, KVR}. En effet l'isomorphisme de Rouvière signifie que
l'application exponentielle réalise un isomorphisme pour la
convolution entre les distributions $K$-invariantes à support $0$
dans $\p$ et les distributions $K$-invariantes à support dans $e\in
G/K$. On peut penser que pour la classe d'espaces symétriques
considérée, cette extension est encore vraie ce qui assurerait que l'isomorphisme de Rouvi\`ere et
notre homomorphisme co\" \i ncident.\\

\end{document}